\documentclass[reqno,11pt]{article}
\usepackage{mathrsfs}
\usepackage{amssymb}
\usepackage{amsthm}
\usepackage{amsmath}
\usepackage{amsfonts}
\usepackage{enumerate}
\usepackage{color}
\usepackage{bm}

\usepackage{graphicx}

\numberwithin{equation}{section}

\newtheorem{theorem}{Theorem}[section]
\newtheorem{lemma}[theorem]{Lemma}
\newtheorem{corollary}[theorem]{Corollary}
\newtheorem{proposition}[theorem]{Proposition}

\theoremstyle{definition}
\newtheorem{definition}[theorem]{Definition}
\newtheorem{assumption}[theorem]{Assumption}
\newtheorem{example}[theorem]{Example}

\theoremstyle{remark}
\newtheorem{remark}[theorem]{Remark}

\begin{document}
\title{Steady Compressible Radially Symmetric Flows with Nonzero Angular Velocity in an Annulus}
\author{Shangkun WENG\thanks{School of mathematics and statistics, Wuhan University, Wuhan, Hubei Province, 430072, People's Republic of China. Email: skweng@whu.edu.cn}\and Zhouping XIN\thanks{The Institute of Mathematical Sciences, The Chinese University of Hong Kong, Shatin, NT, Hong Kong. E-mail: zpxin@ims.cuhk.edu.hk}\and Hongwei YUAN\thanks{The Institute of Mathematical Sciences, The Chinese University of Hong Kong, Shatin, NT, Hong Kong. E-mail: hwyuan@link.cuhk.edu.hk}}
\date{}
\maketitle

\def\be{\begin{eqnarray}}
\def\ee{\end{eqnarray}}
\def\ba{\begin{aligned}}
\def\ea{\end{aligned}}
\def\bay{\begin{array}}
\def\eay{\end{array}}
\def\bca{\begin{cases}}
\def\eca{\end{cases}}
\def\p{\partial}
\def\hphi{\hat{\phi}}
\def\bphi{\bar{\phi}}
\def\no{\nonumber}
\def\eps{\epsilon}
\def\de{\delta}
\def\De{\Delta}
\def\om{\omega}
\def\Om{\Omega}
\def\f{\frac}
\def\th{\theta}
\def\vth{\vartheta}
\def\la{\lambda}
\def\lab{\label}
\def\b{\bigg}
\def\var{\varphi}
\def\na{\nabla}
\def\ka{\kappa}
\def\al{\alpha}
\def\La{\Lambda}
\def\ga{\gamma}
\def\Ga{\Gamma}
\def\ti{\tilde}
\def\wti{\widetilde}
\def\wh{\widehat}
\def\ol{\overline}
\def\ul{\underline}
\def\Th{\Theta}
\def\si{\sigma}
\def\Si{\Sigma}
\def\oo{\infty}
\def\q{\quad}
\def\z{\zeta}
\def\co{\coloneqq}
\def\eqq{\eqqcolon}
\def\di{\displaystyle}
\def\bt{\begin{theorem}}
\def\et{\end{theorem}}
\def\bc{\begin{corollary}}
\def\ec{\end{corollary}}
\def\bl{\begin{lemma}}
\def\el{\end{lemma}}
\def\bp{\begin{proposition}}
\def\ep{\end{proposition}}
\def\br{\begin{remark}}
\def\er{\end{remark}}
\def\bd{\begin{definition}}
\def\ed{\end{definition}}
\def\bpf{\begin{proof}}
\def\epf{\end{proof}}
\def\bex{\begin{example}}
\def\eex{\end{example}}
\def\bq{\begin{question}}
\def\eq{\end{question}}
\def\bas{\begin{assumption}}
\def\eas{\end{assumption}}
\def\ber{\begin{exercise}}
\def\eer{\end{exercise}}
\def\mb{\mathbb}
\def\mbR{\mb{R}}
\def\mbZ{\mb{Z}}
\def\mc{\mathcal}
\def\mcS{\mc{S}}
\def\ms{\mathscr}
\def\lan{\langle}
\def\ran{\rangle}
\def\lb{\llbracket}
\def\rb{\rrbracket}
\def\fr#1#2{{\frac{#1}{#2}}}
\def\dfr#1#2{{\dfrac{#1}{#2}}}
\def\u{{\textbf u}}
\def\v{{\textbf v}}
\def\w{{\textbf w}}
\def\d{{\textbf d}}
\def\nn{{\textbf n}}
\def\x{{\textbf x}}
\def\e{{\textbf e}}
\def\D{{\textbf D}}
\def\U{{\textbf U}}
\def\M{{\textbf M}}
\def\F{{\mathcal F}}
\def\I{{\mathcal I}}
\def\W{{\mathcal W}}
\def\div{{\rm div\,}}
\def\curl{{\rm curl\,}}
\def\R{{\mathbb R}}
\def\FF{{\textbf F}}
\def\A{{\textbf A}}
\def\R{{\textbf R}}
\def\r{{\textbf r}}

\begin{abstract}
  In this paper, we investigate steady inviscid compressible flows with radial symmetry in an annulus. The major concerns are transonic flows with or without shocks. One of the main motivations is to elucidate the role played by the angular velocity in the structure of steady inviscid compressible flows. We give a complete classification of flow patterns in terms of boundary conditions at the inner and outer circle. Due to the nonzero angular velocity, many new flow patterns will appear. There exists accelerating or decelerating smooth transonic flows in an annulus satisfying one side boundary conditions at the inner or outer circle with all sonic points being nonexceptional and noncharacteristically degenerate. More importantly, it is found that besides the well-known supersonic-subsonic shock in a divergent nozzle as in the case without angular velocity, there exists a supersonic-supersonic shock solution, where the downstream state may change smoothly from supersonic to subsonic. Furthermore, there exists a supersonic-sonic shock solution where the shock circle and the sonic circle coincide, which is new and interesting.


  %

%

\end{abstract}

\begin{center}
\begin{minipage}{5.5in}
\end{minipage}
\end{center}

\section{Introduction}\noindent

In this paper, we consider two-dimensional steady compressible Euler flows in an annulus $\Omega=\{(x_1,x_2): r_0<r=\sqrt{x_1^2+x_2^2}<r_1\}$, which are governed by the following system
\begin{align}\label{comeuler2d}
\begin{cases}
\p_{x_1}(\rho u_1)+\p_{x_2}(\rho u_2)=0,\\
\p_{x_1}(\rho u_1^2)+\p_{x_2}(\rho u_1 u_2)+\p_{x_1} p=0,\\
\p_{x_1}(\rho u_1u_2)+\p_{x_2}(\rho u_2^2)+\p_{x_2} p=0,\\
\p_{x_1}(\rho u_1 E+u_1 p)+\p_{x_2}(\rho u_2 E+u_2 p)=0,\\
\end{cases}
\end{align}
where $\u=(u_1,u_2)^t$ is the velocity, $\rho$ is the density, $p$ is the pressure, $E$ is the energy. Here we consider only the polytropic gas, therefore $p= A(S) \rho^{\gamma}$, where $A$ is a smooth function of the entropy $A(S)= R e^{S}$ with $R, \gamma\in (1,+\infty)$ are positive constants. Denote the Bernoulli function by $B=\frac12|{\bf U}|^2+\frac{\ga p}{(\ga-1)\rho}$.

The study on the steady compressible Euler system is not only of fundamental importance in developing the mathematical theory of partial differential equations arising from fluid dynamics, but also has been making a great contribution to the design of projectiles, rockets, aircrafts etc. Courant and Friedrichs \cite[Section 104]{Courant1948} used the Hodograph transformation to rewrite the Euler system \eqref{comeuler2d} into a linear second order PDE on the plane of the flow speed and angle, and obtained some special flows without concerning the boundary conditions. These special flows include circulatory flows and purely radial flows, which are symmetric flows with only angular and radial velocity respectively. Their superpositions are called spiral flows. It had been proved in \cite[Section 104]{Courant1948} that spiral flows can take place only outside a limiting circle where the Jacobian of the hodograph transformation is zero and may change smoothly from subsonic to supersonic or vice verse. However, few results are available concerning properties of transonic spiral flows on the physical plane. This motivates us strongly to investigation of the effects of the angular velocity on general radially symmetric transonic flows with/or without shocks on the physical plane in this paper. The focus will be to find some new transonic flow patterns in the presence of shocks compared with the case without angular velocity.



We will study radially symmetric transonic spiral flows in an annulus with suitable boundary conditions on the inner and outer circle. We will classify all possible transonic radially symmetric flow patterns in terms of physical boundary conditions and study their detailed properties. To this end, we introduce the polar coordinate $(r, \theta)$ as
\begin{eqnarray}
x=r\cos\theta,\ y=r\sin\theta,
\end{eqnarray}
and decompose the velocity as ${\bf u}= U_1 {\bf e}_r + U_2 {\bf e}_{\theta}$ with
\begin{align}
\e_r=\begin{pmatrix} \cos\theta\\ \sin\th\end{pmatrix},
\e_\th=\begin{pmatrix}-\sin\th\\ \cos\th\end{pmatrix}.
\end{align}
Therefore the system \eqref{comeuler2d} can be rewritten as
\begin{eqnarray}\label{ProblemImm}
\begin{cases}
\partial_r(\rho U_1)+\frac{1}{r}\partial_{\theta}(\rho U_2)+\frac{1}{r} \rho U_1=0,\\
(U_1\partial_r +\frac{U_2}{r}\partial_{\theta}) U_1+\frac{1}{\rho} \partial_r p-\frac{U_2^2}r=0,\\
(U_1\partial_r +\frac{U_2}{r}\partial_{\theta}) U_2+\frac{1}{r\rho} \partial_{\theta}p+\frac{U_1U_2}r=0,\\
(U_1\partial_r +\frac{U_2}{r}\partial_{\theta}) A=0.
\end{cases}
\end{eqnarray}

For the radially symmetric solutions to \eqref{ProblemImm} of the form ${\bf u}=U_{1}(r){\bf e}_r+U_{2}(r){\bf e}_\theta$, $\rho=\rho(r)$, $A=A(r)$ in $\Omega$, the steady Euler system \eqref{ProblemImm} reduces to
\begin{align}\label{radial-euler}
\begin{cases}
\frac{d}{dr}(\rho U_1)+\frac1r \rho U_1=0, \ \ &r_0<r<r_1,\\
U_1U_1'+\frac1\rho \frac{d}{dr}p-\frac{U_2^2}r=0,\ \ &r_0<r<r_1,\\
U_1U_2'+\frac{U_1U_2}r=0,\ \ &r_0<r<r_1,\\
U_1A'=0,\ \ &r_0<r<r_1.
\end{cases}
\end{align}
We start with smooth solutions to \eqref{radial-euler} with one side boundary conditions prescribed on the inner or outer circle and give a classification of possible flow patterns in an annulus including purely smooth supersonic, subsonic flows, and smooth decelerating or accelerating transonic flows. We further analyze the detailed properties of these solutions and their dependence on the boundary data. If the boundary data is prescribed on the outer circle and the fluid moves from outer circle to the inner one, there exists a limiting circle such that the acceleration will blow up when the fluid moves to the limiting circle. The sonic circle and the limiting circle are also determined by boundary data. Although these smooth transonic spiral flows are essentially same as the ones obtained by Courant-Friedrichs \cite{Courant1948} on the hodograph plane, the detailed behaviors of these solutions would be very helpful for the investigation of structural stability of symmetric transonic flows under suitable perturbations on the boundaries \cite{wxy20}.

An important fact is on the degeneracy type of the sonic points in the smooth radially symmetric transonic spiral flows. Recall the definition by Bers \cite{Bers1958}, a sonic point in a $C^2$ smooth transonic flow is exceptional if and only if the velocity is orthogonal to the sonic curve at this point. Then all sonic points of smooth transonic spiral flows here are nonexceptional and noncharacteristically degenerate due to the nonzero angular velocity. This is different from smooth transonic flows constructed by Wang and Xin \cite{WX2013,WX2016,Wang2019d,Wang2020} in De Laval nozzles. In particular, the authors in \cite{Wang2019d} proved that if the nozzle is suitably flat at its throat, then there exists a unique smooth transonic irrotational flow of Meyer type with all sonic points being exceptional in De Laval nozzles. The sonic points must be located at the throat of the nozzle and are strongly degenerate in the sense that all the characteristics from sonic points coincide with the sonic line and can not approach the supersonic region. It should be mentioned that Kuzmin \cite{Kuzmin2002} had studied the structural stability of accelerating smooth transonic flows with some artificial boundary conditions on the potential and stream function plane. The existence of subsonic-sonic weak solutions to the 2-D steady potential equation were proved in \cite{cdsw07,xx07} by utilizing the compensated compactness and later on the subsonic-sonic limit for multidimensional potential flows and steady Euler flows were examined in \cite{chw16,hww11}. However, the solutions obtained by the subsonic-sonic limit only satisfied the equations in the sense of distribution and the regularity of the solutions and the properties of sonic points are not clear at all.


More importantly, we study transonic solutions with shocks to \eqref{ProblemImm} by prescribing suitable boundary conditions on the inner and out circle. Recall that, a piecewise smooth solution $(U_{1}^{\pm},U_{2}^{\pm}, \rho^{\pm}, A^{\pm})\in C^1(\Omega^{\pm})$ with a jump on the curve $r=f(\th)$ is a shock solution to the system \eqref{ProblemImm}, if $(U_{1}^{\pm},U_{2}^{\pm}, \rho^{\pm}, A^{\pm})$ satisfy the system \eqref{ProblemImm} in $\Omega^{\pm}$ respectively, and the entropy condition $[p]>0$ and the following Rankine-Hugoniot jump conditions hold on $r=f(\theta)$:
\begin{align}\label{jump00}
\begin{cases}
[\rho U_1]-\frac{f'(\theta)}{f(\theta)}[\rho U_2]=0,\\
[\rho U_1^2+p]-\frac{f'(\theta)}{f(\theta)}[\rho U_1U_2]=0,\\
[\rho U_1U_2]-\frac{f'(\theta)}{f(\theta)}[\rho U_2^2+p]=0,\\
[B]=0,\\
\end{cases}
\end{align}
where $[v]=v^+(f(\theta),\theta)-v^-(f(\theta),\theta)$. For radially symmetric transonic solutions with shocks, we would consider the system \eqref{radial-euler} supplemented with suitable boundary conditions on both inner and outer circle and the shock curve will be a circle $r=r_b$ and the Rankine-Hugoniot jump conditions reduce to
\begin{align}\label{rankine-hugoniotII}
\begin{cases}
[\rho U_1](r_b)=[\rho U_1^2+p](r_b)=0,\\
[\rho U_1U_2](r_b)=[B](r_b)=0.
\end{cases}
\end{align}


The existence and uniqueness of radially symmetric supersonic-subsonic shock solutions without angular velocity in a divergent nozzle have been proved in \cite{Courant1948,Xin2008a}. The dynamical and structural stability of these symmetric transonic shock solutions have been an important and difficult research topic in the mathematical studies of gas dynamics. Symmetric transonic shocks are shown in \cite{Xin2008a} to be dynamically stable in divergent nozzles and is dynamically unstable in convergent nozzles. The transonic shock problem with given exit pressure was shown to be ill-posed for the potential flow model in \cite{xy05,xy08a}. The authors in \cite{lxy09b,LiXin2013} proved that, for the two-dimensional steady compressible Euler system, the symmetric transonic shock solutions in \cite{Courant1948,Xin2008a} are structurally stable under the perturbation of the exit pressure and the nozzle wall. They further investigated the existence and monotonicity of the axi-symmetric transonic shock by perturbing axi-symmetrically exit pressure in \cite{lxy10a,lxy10b}. The authors in \cite{Weng2019} examined the structural stability under axi-symmetric perturbations of the nozzle wall by introducing a modified invertible Lagrangian transformation. There have been many interesting results on transonic shock problems in flat or divergent nozzles for different models with various exit boundary conditions, such as the non-isentropic potential model, the exit boundary condition for the normal velocity, the spherical symmetric flows without solid boundary, etc. One may refer to \cite{Bae2011,cf03,ccs06,ccf07,chen05,chen08,cy08,FangGao,FangXin,lxy09a,Liu2016} and the references therein.

However, in the presence of the nonzero angular velocity, many new wave patterns appear in the annulus. It follows from \eqref{rankine-hugoniotII} that the radial velocity will jump from supersonic to subsonic and the angular velocity experiences no jump across the shock. Therefore the total velocity after the shock may be supersonic, sonic and subsonic. We will give a complete classification of flow patterns according to the exit pressure and the position of the outer circle. Besides the well-known supersonic-subsonic shock solutions in a divergent nozzle as in the case without the angular velocity, we find that there are supersonic-supersonic shock solutions and the flow at downstream may change smoothly from supersonic to subsonic. Furthermore, there exists also a supersonic-sonic shock solution where the shock front and the sonic circle coincide, which is new and surprising.

The structure of the paper will be organized as follows. In Section \ref{smooth}, we will consider the one side boundary value problem to \eqref{radial-euler} and show that there are accelerating or decelerating smooth transonic flows. In Section \ref{transonic}, the precise description of transonic flows with shocks and new wave patterns will be presented by studying \eqref{radial-euler} with two sides boundary values.

\section{Smooth radially symmetric flows with nonzero angular velocity}\label{smooth}\noindent

In this section, we will construct a radially symmetric smooth solution to \eqref{comeuler2d} in an annulus $\Omega=\{(x_1,x_2): r_0<r=\sqrt{x_1^2+x_2^2}<r_1\}$, where the fluid flows from the inner circle to the outer one. Thus we consider the following problem:

{\bf Problem I.} Find smooth functions ${\bf u}(x)= U_{1}(r) {\bf e}_r + U_{2}(r) {\bf e}_{\theta}, \rho(x)=\rho(r)$ and $p(x)=p(r), A(x)=A(r)$, which solve the system \eqref{radial-euler} with the boundary conditions on $r=r_0$:
\begin{eqnarray}\label{ProblemI}
\rho(r_0)=\rho_{0},\ U_{1}(r_0)=U_{10}>0, \ U_{2}(r_0)=U_{20},\ A(r_0)=A_{0}.
\end{eqnarray}

It is easy to see that the problem \eqref{radial-euler} with \eqref{ProblemI} is equivalent to
\begin{align}\label{alg1}
\begin{cases}
\rho U_1=\dfr{\ka_1}{r},\ \kappa_1= r_0 \rho_{0} U_{10},\\
B\equiv B_0, B_0=\dfrac{1}{2}(U_{10}^2+U_{20}^2)+\dfr{A_0 \gamma}{\gamma-1} \rho_{0}^{\gamma-1},\\
U_2=\dfr{\ka_2}{r},\ \kappa_2= r_0 U_{20},\\
A\equiv A_0.
\end{cases}
\end{align}

\bt\label{ProbI}
There are three distinguished cases.
\begin{enumerate}[(1)]
  \item If the incoming flow at the entrance is supersonic in the $r$-direction, that is $U_1^2(r_0)>c^2(r_0)$, then there exists a unique supersonic smooth solution to Problem I in $\Omega$.
  \item If the incoming flow at the entrance is subsonic, that is $|{\bf U}|^2(r_0)<c^2(r_0)$, then there exists a unique subsonic smooth solution to Problem I in $\Omega$.
  \item If the incoming flow at the entrance is supersonic but subsonic in the $r$-direction, that is $|{\bf U}|^2(r_0)>c^2(r_0)>U_1^2(r_0)$, there exists a constant $r_c>r_0$ defined in \eqref{critical-r}, which depends only on $r_0,\gamma$ and the incoming flow, such that
  \begin{enumerate}
  \item if $r_1>r_c$, there exists a unique smooth transonic solution to Problem I in $\Omega$;
  \item if $r_1<r_c$, there exists a unique smooth supersonic solution to Problem I in $\Omega$;
  \item if $r_1=r_c$, there exists a unique smooth supersonic-sonic solution to Problem I in $\Omega$ with sonic circle located at $r=r_c$.
  \end{enumerate}
\end{enumerate}
\et

\begin{remark}
{\it If $c^2(r_0)=U_1^2(r_0)$, then the derivatives of $\rho$ and $U_1$ at $r_0$ become infinite. If $|{\bf U}|^2(r_0)=c^2(r_0)>U_1^2(r_0)$, then there exists a unique smooth sonic-subsonic solution to Problem I in $\Omega$ with sonic circle located at $r=r_0$.}
\end{remark}

\begin{remark}
{\it If $U_{10}=0$, then $U_1\equiv 0, A(r)\equiv A_0$, $U_2(r)=\frac{\kappa_2}{r}$ and
\begin{eqnarray}\nonumber
\rho(r)=\left(\frac{\gamma-1}{A_0\gamma}\right)^{\frac{1}{\gamma-1}} \left(B_0-\frac{\kappa_2^2}{2r^2}\right)^{\frac{1}{\gamma-1}}.
\end{eqnarray}
If $U_{20}^2< A_0\gamma \rho_0^{\gamma-1}$, then the flow is purely subsonic in $[r_0,\infty)$. If $U_{20}^2>A_0\gamma \rho_0^{\gamma-1}$, then the flow changes smoothly from supersonic to subsonic as $r>r_0$ goes to infinity with a sonic circle located at $r=r_c=\sqrt{\frac{\gamma+1}{2(\gamma-1) B_0}}|\kappa_2|$.
}\end{remark}

\br
{\it
In \cite{WX2016}, the authors have analyzed the geometric property and degeneracy types of sonic points for general $C^2$ smooth subsonic-sonic flows and transonic steady potential flows. It is shown that a sonic point satisfying the interior subsonic circle condition is exceptional if and only if the governing equation is characteristically degenerate at this point. Moreover, they proved the existence and uniqueness of smooth transonic flows of Meyer type in \cite{Wang2019d}, in which all sonic points are exceptional and characteristically degenerate. It should be remarked here that the smooth transonic flows constructed in Theorem \ref{ProbI} have nonzero angular velocities, and thus all sonic points are nonexceptional and noncharacteristically degenerate. The structural stability of this special solution will be investigated in \cite{wxy20}.
}\er

Denote $M_i^2:=\frac{U_i^2}{c^2}$ for $i=1,2$, ${\bf M}=(M_1,M_2)^t$ and set
\be\no
K(r)=\frac{2(\gamma-1)}{\ga+1}B-\fr{\ga-1}{\ga+1}U_2^2(r).
\ee
Now we prove Theorem \ref{ProbI}.

\bpf

Direct calculations show that
\be\label{Mp0}
\rho'&=&\dfr{|\M|^2}{r(1-M_{1}^2)}\rho ,\ U_{1}'=-\dfr{1+ M_{2}^2}{r(1-M_{1}^2)}U_{1},\ \ U_2'=-\dfr{U_2}r,\\\label{Mp1}
(M_{1}^2)'&=&\dfr{M_{1}^2}{r(M_{1}^2-1)}\bigg(2+(\ga-1)M_{1}^2+(\ga+1)M_{2}^2\bigg)\\\label{Mp2}
(M_{2}^2)'&=&\dfr{M_{2}^2}{r(M_{1}^2-1)}\bigg(2+(\ga-3)M_{1}^2+(\ga-1)M_{2}^2\bigg)\\\label{Mpt}
(|\M|^2)'&=&\dfr{|\M|^2}{r(M_{1}^2-1)}\bigg(2+(\ga-1)|\M|^2\bigg).
\ee
\begin{enumerate}[(1)]
  \item Suppose that $M_1^2(r_0)-1>0$. Then by \eqref{Mp1}, $M_1^2(r)-1$ increases as $r$ increases. So $|{\bf M}|^2(r)>M_1^2(r)>1$ for any $r>r_0$ and the flow is always supersonic. The existence and uniqueness of the solution follow directly from the standard theory of ODE system. Furthermore, one has $\frac{d}{dr}\log \rho<-\frac{1}{r}$, which implies that $\rho(r)\leq \rho_0\dfr{r_0}r\rightarrow 0$ as $r\to \infty$. Also $|U_2(r)|,|U'_2(r)|\rightarrow0$ as $r\rightarrow +\infty$, by the Bernoulli's law, $U_1(r)\rightarrow\sqrt{2B_0}$. \eqref{Mp1} implies that $(M_{1}^2)'\geq \dfr{\ga-1}r M_1^2$, which yields
  \be\label{M1infty}
  M_1^2(r)\geq M_1^2(r_0) (\dfr{r}{r_0})^{\ga-1}.
  \ee

  \item Suppose that $M_1^2(r_0)<|\M|^2(r_0)<1$. Then by \eqref{Mp1}, $\frac{d}{dr} M_1^2(r)<0$ and thus $M_1^2(r)<1$ for any $r>r_0$. This together with \eqref{Mpt} and $|\M|^2(r_0)<1$, implies that $|\M|^2(r)<1$ for any $r>r_0$ and the flow is always subsonic. Furthermore, $|{\bf M}|^2(r)\leq \frac{r_0^2|{\bf M}|^2(r_0)}{r^2}\rightarrow0$, $U_1(r)\rightarrow0$, $c^2(r)\rightarrow(\ga-1)B_0$, $|U_2(r)|,|U'_2(r)|\rightarrow0$ as $r\rightarrow +\infty$.

  \item Suppose that $M_1^2(r_0)<1<|\M|^2(r_0)$. Then one can get a unique smooth solution with $M_1^2(r)<1$ for any $r>r_0$ by same arguments as in (2). Also there holds $|{\bf M}|^2(r)\leq \frac{r_0^2|{\bf M}|^2(r_0)}{r^2}\rightarrow0$ as $r\to \infty$. The position of the sonic circle can be determined as follows. Define the critical density
  \begin{align}\label{rhostar}
  \rho_c=\bigg(\frac{2(\ga-1)B_0}{(\ga+1)\ga A_0}\bigg)^{\fr1{\ga-1}}.
  \end{align}
  Note that $U_2$ and $A$ are already obtained. Since $U_1=\frac{\kappa_1}{r \rho}$, it is easy to see that $\rho$ is a root of the following algebraic equation
  \begin{align}\label{Fr}
  F_r(\rho):=\dfr{\ga}{\ga-1}A_0\rho^{\ga+1}-(B_0-\dfr12
  \dfr{\kappa_2^2}{r^2})\rho^2+\dfr{\kappa_1^2}{2r^2}=0.
  \end{align}
  The sonic circle $r=r_c$ is a root of $F_{r}(\rho_c)=0$, which can be calculated explicitly
  \begin{eqnarray}\label{critical-r}
  r_c=\sqrt{\frac{(\gamma+1)(\kappa_1^2+ \kappa_2^2\rho_c^2)}{2(\gamma-1)B_0\rho_c^2}}.
  \end{eqnarray}
  Therefore if $r_1>r_c$ one obtains a transonic smooth solution on $\Omega$; if $r_1=r_c$ one gets a supersonic-sonic flow; if $r_1<r_c$, there exists a purely hyperbolic flow.

\end{enumerate}

\epf

Next, we consider the existence of radially symmetric smooth solutions to \eqref{comeuler2d} in an annulus $\Omega$, where the fluid moves from the outer circle to the inner circle. That is

{\bf Problem II.} Find smooth functions ${\bf u}(x)= U_{1}(r) {\bf e}_r + U_{2}(r) {\bf e}_{\theta}, \rho(x)=\rho(r)$ and $p(x)=p(r), A(x)=A(r)$, which solve the system \eqref{radial-euler} supplemented with the boundary conditions at the outer circle
\begin{eqnarray}\label{ProblemI2}
\rho(r_1)=\rho_{0},\ U_{1}(r_1)=U_{10}<0, \ U_{2}(r_1)=U_{20},\ A(r_1)=A_{0}.
\end{eqnarray}

Slightly different from {\bf Problem I}, there exists a limiting circle $r=r^{\sharp}<r_1$ such that the acceleration of the solution to {\bf Problem II} will blow up as $r$ tends to $r^{\sharp}$. The main results for {\bf Problem II} can be stated as follows.

\bt\label{ProbII}
{\it For a given incoming flow at the outer circle $r=r_1$, there exists a limiting circle $r^\sharp<r_1$ depending only on $r_1, \gamma$ and the incoming flow state at $r_1$ such that the flow develops singularities when $r\to r^{\sharp}$ in the sense that if $r_0<r^\sharp$ there does not exist any global $C^1$-smooth solution to Problem II in $\Omega$. If $r_0\geq r^\sharp$, there exists a smooth solution to Problem II in $\Omega$, the flow pattern in $\Omega$ can be divided into the following cases.
\begin{enumerate}[(1)]
  \item Suppose that $M_1^2(r_1)>1$ or $|{\bf M}|^2(r_1)>1>M_1^2(r_1)$, then there exists an unique smooth supersonic solution to Problem II in $\Omega$.
  \item Suppose that $|{\bf M}|^2(r_1)<1$, there exists a constant $0<r=r_c<r_1$ defined in \eqref{critical-r}, which depends only $r_1,\gamma$ and the incoming flow, such that
      \begin{enumerate}[(a)]
      \item if $r_0>r_c$, there exists a unique smooth subsonic solution to Problem II in $\Omega$;
      \item if $r_0<r_c$, there exists a unique smooth transonic solution to Problem II in $\Omega$;
      \item if $r_0=r_c$, there exists an unique smooth subsonic-sonic solution to Problem II in $\Omega$ with sonic circle located at $r=r_c$.
      \end{enumerate}
\end{enumerate}
}\et

\br
{\it If $c^2(r_1)=U_1^2(r_1)$, then the derivatives of $\rho$ and $U_1$ at $r_1$ become infinite. If $|{\bf U}|^2(r_1)=c^2(r_1)>U_1^2(r_1)$, then there exists a unique smooth sonic-supersonic solution to Problem II in $\Omega$ with the sonic circle located at $r=r_1$.}
\er

\br\no
{\it The existence results of smooth supersonic and transonic flows in Theorem \ref{ProbII} are essentially same as the ones obtained by Courant-Friedrichs \cite{Courant1948} on the hodograph plane. Here we study the boundary value problem \eqref{radial-euler} with one side boundary conditions on the physical plane, and analyze the dependence of the solutions on boundary data, which also suits our purpose for the investigation of structural stability of symmetric transonic flows under suitable perturbations of boundary conditions \cite{wxy20}.
}\er

Note that for smooth solutions,  {\bf Problem II} and \eqref{alg1} are equivalent, where $\kappa_1= r_1 \rho_{10} U_{10}$ and $\kappa_2= r_1 U_{20}$.
\bl\label{Lemma2}
There exists a $0<r^\sharp< r_1$ such that
\begin{enumerate}[(1)]
  \item for any $r\in(r^\sharp,r_1)$, the algebraic system \eqref{alg1} has exactly two solutions $(\rho^-(r),U_1^-(r),U_2(r),A(r))$ and $(\rho^+(r),U_1^+(r),U_2(r),A(r))$ with $(M_1^-)^2(r)>1$ and $(M_1^+)^2(r)<1$.
  \item for $r=r^\sharp$, the algebraic system \eqref{alg1} has exactly one solution $(\rho^-(r^\sharp),U_1^-(r^\sharp),U_2(r^\sharp),A(r^\sharp))=(\rho^+(r^\sharp),U_1^+(r^\sharp),U_2(r^\sharp),A(r^\sharp))$ with $(M_1^-)^2(r^\sharp)=1=(M_1^+)^2(r^\sharp)$.
  \item for any $r\in(0,r^\sharp)$, the algebraic system \eqref{alg1} has no solution.
\end{enumerate}
\el
\bpf
Recall that $\rho$ satisfies the algebraic equation $F_r(\rho)=0$ with $F_r(\rho)$ defined by \eqref{Fr}. Set $\tilde{r}=r_1 U_{20}/\sqrt{2B_0}$. Then for each $0<r\leq \tilde{r}$, $F_r(\rho)$ is monotonically increasing in $\rho$, $F_r(\rho)\geq F_r(0)=\ka_1^2/(2r^2)>0$ and $F_r(\rho)=0$ has no solution in $[0,+\infty)$.

For each $r>\tilde{r}$, one has
\begin{align}
F_r'(\rho)=\fr{\ga(\ga+1)}{\ga-1}A_0\rho^{\ga}-2(B_0-\fr12 \fr{\kappa_2^2}{r^2})\rho.
\end{align}
Then $F_r(\rho)$ attains its minimum at $\rho=\rho_*(r)$, where
\begin{align}\label{rhostard}
\rho_*(r)=\bigg(\fr{K(r)}{\ga A_0}\bigg)^{\fr1{\ga-1}}.
\end{align}

Some direct computations yield
\begin{align*}
&F_r(\rho_*(r))=-\fr12(\ga A_0)^{-\fr{2}{\ga-1}}K(r)^{\fr{\ga+1}{\ga-1}}+\fr12r^{-2}\kappa_1^2,\\
&F_{\tilde{r}}(\rho_*(\tilde{r}))=\frac{\kappa_1^2}{2\tilde{r}^{2}}>0,\\
&F_{r_1}(\rho_*(r_{1}))=-\fr12\rho_0^2U_{10}^2\left(\left(\fr{\ga-1}{\ga+1}M_{10}^{\fr4{\ga+1}}+\fr2{\ga+1}M_{10}^{-\fr{2(\ga-1)}{\ga+1}}\right)^{\fr{\ga+1}{\ga-1}}-1\right),\\
&\fr{d}{dr}F_r(\rho_*(r))=-\left((\ga A_0)^{-\fr2{\ga-1}}K(r)^{\fr2{\ga-1}}\kappa_2^2+\kappa_1^2\right)r^{-3}<0, \ \ \forall r>\tilde{r}.
\end{align*}
Thus $F_r(\rho_*(r))$ is monotonically decreasing as a function of $r$. Furthermore, Young's inequality implies that $F_{r_1}(\rho_*(r_{1}))<0$ whenever $M_{10}^2\neq 1$. Hence there exists a unique $r^\sharp\in (\tilde{r},r_1)$ such that $F_{r^\sharp}(\rho_*(r^{\sharp}))=0$. And if $r\in(\tilde{r},r^\sharp)$, $F_r(\rho_*(r))>0$; if $r\in(r^\sharp,r_1)$, then $F_r(\rho_*(r))<0$. Note that for each $r\in [r^{\sharp}, r_1)$, $F_r(\rho)$ decreases first as $\rho$ increases from $0$ to $\rho_*(r)$ and then increases when $\rho$ increases further. Therefore for each $r\in(r^\sharp,r_1)$, $F_r(\rho)=0$ has exactly two solution $\rho^\pm(r)$ with $0<\rho^-(r)<\rho_*(r)<\rho^+(r)$. For $r=r^\sharp$, $F_r(\rho)=0$ has exactly one solution $\rho=\rho_*(r)$. For each $r\in (\tilde{r}, r^{\sharp})$, $F_r(\rho)=0$ has no solution in $[0,+\infty)$. In addition, $\rho^-(r)<\rho_*(r)<\rho^+(r)$ is equivalent to $M_1^+(r)<1<M_1^-(r)$.

\epf

By Lemma \ref{Lemma2} and a similar argument as in the proof of Theorem \ref{ProbI}, one can establish Theorem \ref{ProbII}.

\br\no
{\it Note that $\tilde{r}$ is the position where $c^2(\tilde{r})=U_1^2(\tilde{r})=0$ and $U_2^2(\tilde{r})=2B_0$.
}\er

\br\label{estirsharp}
{\it $r^\sharp$ is the position where $U_1^2(r^\sharp)=c^2(r^\sharp)$ and there is another estimate of $r^\sharp$ as follows. $U_1^2(r^\sharp)=c^2(r^\sharp)$ implies $\ga A_0 (\fr{\kappa_1}{r^\sharp})^{\ga-1}=U_1^{\ga+1}(r^\sharp)=(\fr{2(\ga-1)B_0}{\ga+1}-\fr{\ga-1}{\ga+1}U_2^2(r^\sharp))^{(\ga+1)/2}
\leq(\fr{2(\ga-1)B_0}{\ga+1})^{(\ga+1)/2}$. Thus $r^\sharp \geq (\ga A_0)^{1/(\ga-1)}(\fr{\ga+1}{2(\ga-1)B_0})^{\fr{\ga+1}{2(\ga-1)}}\kappa_1$.
}
\er

\section{Transonic shock flows in an annulus}\label{transonic}\noindent

In this section, we will prove the existence of radially symmetric transonic solutions with shocks and nonzero angular velocity to \eqref{radial-euler} satisfying suitable boundary conditions both at the inner and outer circle. Especially, motivated by previous studies due to Courant-Friedrichs\cite{Courant1948}, we prescribe the pressure at the exit. Therefore the problem to be solved can be formulated precisely as

{\bf Problem III.} Construct a piecewise smooth radially symmetric flow on an annulus $\Omega=\{(x_1,x_2):r_0<r=\sqrt{x_1^2+x_2^2}<r_1\}$ which moves from the inner circle to the outer one with a shock at $r=r_b\in (r_0,r_1)$. More precisely, one looks for smooth functions $(U_{1}^{\pm},U_{2}^{\pm}, \rho^{\pm}, p^{\pm}, A^{\pm}$ on $\Omega^{\pm}$ respectively with $\Omega^+=\{(x_1,x_2): r_b<r<r_1\}$ and $\Omega^-=\{(x_1,x_2): r_0<r<r_b\}$, which solves \eqref{radial-euler} on $\Omega^{\pm}$ with a discontinuity $r=r_b$ satisfying the physical entropy condition $[p]>0$ and the following Rankine-Hugoniot conditions \eqref{rankine-hugoniotII} and satisfy the boundary conditions
\begin{eqnarray}\label{bd1}
&&\rho(r_0)=\rho_{0},\ U_{1}(r_0)=U_{10}>0, \ U_{2}(r_0)=U_{20}(\neq 0),\ A(r_0)=A_{0},\\\label{bd2}
&&p(r_1)= p_{ex}.
\end{eqnarray}

By Theorem \ref{ProbI}, there exists a unique smooth solution $(U_1^-, U_2^-, \rho^-, A^-)$ to \eqref{radial-euler} with the boundary condition \eqref{bd1}. Define
\be\label{f1}
f_1(r)&=&\bigg(1-\dfr{(U_2^-)^2}{\fr{2(\ga-1)}{\ga+1}B_0-\fr{\ga-1}{\ga+1}(U_2^-)^2}\bigg)\dfr{(U_1^-)^2}{(c^-)^2},\\\label{f2}
f_2(\rho)&=&\dfr{\ga p_{ex}}{\ga-1}\rho-(B_0-\dfr12\dfr{\ka_2^2}{r_1^2})\rho^2+\dfr{\ka_1^2}{2r_1^2},
\ee
and
\be\nonumber
r_*=\sqrt{\frac{\ga \kappa_2^2}{(\ga-1)B_0}},\ \ \rho^{\sharp}=\dfrac{\ga(\ga+1) p_{ex}}{2(\ga-1)B_0}.
\ee
Suppose that $\kappa_2\neq 0$. It will be proved that $f_1(r)$ increases monotonically in $r>r_*$. Since $f_1(r_*)=0$ and $\displaystyle\lim_{r\to \infty} f_1(r)=\infty$, there exists a unique $r'_*>r_*$ such that $f_1(r'_*)=1$. Note that $r'_*$ is uniquely determined by the incoming flow and the system \eqref{radial-euler}.

Now the main results can be stated as follows.

\bt\label{ProbIII}
{\it For a given incoming flow which is supersonic in the $r$-direction and has a nonzero angular velocity at the entrance, there exist two constants $p_1<p_0$ such that, for any $p_{ex}\in(p_1,p_0)$, there exist a unique piecewise smooth solution to Problem III in $\Omega$ with a shock located at $r\equiv r_b$ and satisfying $U_1^-(r)>c^-(r)$ for $r\in(r_0,r_b)$, $U_1^+(r)<c^+(r)$ for $r\in(r_b,r_1)$. Furthermore, the solution has the following properties.
\begin{enumerate}[(A)]
  \item The shock position $r=r_b$ increases as the exit pressure $p_{ex}$ decreases. In addition, the shock position $r_b$ tends to $r_1$ if $p_{ex}$ goes to $p_1$ while $r_b$ approaches to $r_0$ if $p_{ex}$ goes to $p_0$.
  \item The flow patterns in $\Omega^+$ can be classified in terms of the boundary conditions as follows:
  \begin{enumerate}[1.]
  \item If $r_0<r_1\leq r'_*$ and $p_{ex}\in(p_1,p_0)$, then
  \begin{description}
  \item{Subcase 1.1.} for $f_2(\rho^\sharp)>0$, there exists a supersonic-supersonic shock and the flow is supersonic in $\Omega^+$.
  \item{Subcase 1.2.} for $f_2(\rho^\sharp)<0$, there exists a supersonic-supersonic shock and the flow changes smoothly from supersonic to subsonic in $\Omega^+$.
  \item{Subcase 1.3.} for $f_2(\rho^\sharp)=0$, there exists a supersonic-supersonic shock and the flow is supersonic in $\Omega^+$ but degenerates to a sonic state at the exit.
  \end{description}
  \item If $r_0<r'_*<r_1$, then there exists a $p'_*\in(p_1,p_0)$ such that,
  \begin{description}
  \item{Subcase 2.1.} for $p_{ex}\in(p'_*,p_0)$ with $f_2(\rho^\sharp)>0$, a supersonic-supersonic shock exists and the flow is supersonic in $\Omega^+$.
  \item{Subcase 2.2.} for $p_{ex}\in(p'_*,p_0)$ with $f_2(\rho^\sharp)<0$, a supersonic-supersonic shock exists and the flow changes smoothly from supersonic to subsonic in $\Omega^+$.
  \item{Subcase 2.3.} for $p_{ex}\in(p'_*,p_0)$ with $f_2(\rho^\sharp)=0$, a supersonic-supersonic shock exists and the flow is supersonic in $\Omega^+$ but degenerates to a sonic state at the exit.
  \item{Subcase 2.4.} for $p_{ex}\in(p_1,p'_*)$, a supersonic-subsonic shock exists and the flow is subsonic in $\Omega^+$.
  \item{Subcase 2.5.} for $p_{ex}=p'_*$, there exists a supersonic-sonic shock and the flow is subsonic in $\Omega^+$ but degenerates to a sonic state at the shock position.
  \end{description}
  \item If $ r_1>r_0\geq r'_*$, then for any $p_{ex}\in(p_1,p_0)$, there exists a supersonic-subsonic shock and the flow is uniformly subsonic in $\Omega^+$.
\end{enumerate}
\end{enumerate}
}
\et
\begin{remark}\nonumber
{\it The value $f_2(\rho^\sharp)$ is completely determined by the boundary conditions, while $f_1(r_1)$ is determined by the incoming flow and the steady Euler system a priori.
}\end{remark}

\begin{remark}\nonumber
{\it In subcase 2.5, there is a new flow pattern that is not observed before, where the shock front and the sonic curve coincide due to the nonzero angular velocity. This supersonic-sonic shock is structurally unstable in general.
}\end{remark}

\begin{remark}\nonumber
{\it If $\kappa_2=0$, then $r_*=0$ and $f_1(r)=(M_1^-(r))^2$. Therefore $f_1(r)>1$ for any $r>r_0$ and only the case 3 in Theorem \ref{ProbIII} happens. This coincides with the results obtained in \cite{Xin2008a}.
}
\end{remark}

\begin{remark}\nonumber
{\it The stability of these transonic flows with shocks will be studied in the future.
}
\end{remark}

\subsection{Proof of Theorem \ref{ProbIII}}\noindent


In $\Omega^-$, it holds that
\begin{align}
\begin{cases}
\rho^- U_1^-=\dfr{\ka_1}{r},\ \kappa_1= r_0 \rho_{0} U_{10},\\
B^-\equiv B_0, B_0=\dfrac{1}{2}(U_{10}^2+U_{20}^2)+\dfrac{A_0 \gamma}{\gamma-1} \rho_{0}^{\gamma-1},\\
U_2^-=\dfr{\ka_2}{r},\ \kappa_2= r_0 U_{20},\\
A^-\equiv A_0\equiv A_0^-.
\end{cases}
\end{align}
It follows from \eqref{rankine-hugoniotII} that $\ka_1$, $\ka_2$ and $B_0$ are unchanged across the shock. Therefore, in $\Omega^+$, one can get that
\begin{align}\label{Omgp}
\begin{cases}
(\rho^+ U_1^+)(r)=\dfr{\ka_1}{r},\ \ U_2^+(r)=\dfr{\ka_2}{r},\\
B^+\equiv B_0,\ \ A^+\equiv A_0^+,
\end{cases}
\end{align}
where $A_0^+$ is to be determined by \eqref{rankine-hugoniotII}.

First, we solve \eqref{rankine-hugoniotII} for fixed $r_b\in(r_0,r_1)$. It follows from the first two equations in \eqref{rankine-hugoniotII} that $[U_1+\fr{c^2}{\ga U_1}]=0$, which further implies that $[U_1+\fr{K_0}{U_1}]=0$ with $K_0=K(r_b)$, where $K(r)=\frac{2(\gamma-1)}{\ga+1}B-\fr{\ga-1}{\ga+1}U_2^2(r)$ has no jump across the shock since $\ka_2$ and $B_0$ are unchanged across the shock. Thus on the shock
\be\label{U1p}
U_1^+(r_b)=\dfr{K_0}{U_1^-(r_b)}.
\ee
If $U_1^-(r_b)>c^-(r_b)$, then $U_1^+(r_b)<c^+(r_b)$. By the first equation in \eqref{Omgp},
\be\label{rhop}
\rho^+(r_b)=\dfr{\ka_1U_1^-(r_b)}{r_b K_0}=\dfr{\ka_1^2}{K_0r_b^2\rho^-(r_b)}.
\ee

Moreover, one can conclude from $[B]=0$ and $[\rho U_1^2+ p]=0$ that
\begin{align*}
\begin{cases}
A_0^+(U_1^-(r_b))^{\ga-1}-A_0^-(U_1^+(r_b))^{\ga-1}=(\fr{r_b}{\kappa_1})^{\ga-1}\fr{\ga-1}{2\ga}((U_1^-(r_b))^2-(U_1^+(r_b))^2)(U_1^-(r_b)U_1^+(r_b))^{\ga-1}\\\no
A_0^+(U_1^-(r_b))^{\ga}-A_0^-(U_1^+(r_b))^{\ga}=(\fr{r_b}{\kappa_1})^{\ga-1}(U_1^-(r_b)-U_1^+(r_b))(U_1^-(r_b)U_1^+(r_b))^{\ga}.
\end{cases}
\end{align*}
Thus
\begin{align}\no
A_0^+=&(\fr{r_b}{\kappa_1})^{\ga-1} (U_1^+(r_b))^{\ga}(\fr{\ga+1}{2\ga}U_1^-(r_b)-\fr{\ga-1}{2\ga}U_1^+(r_b))\\\no
=&(\fr{r_b}{\kappa_1})^{\ga-1} (\dfr{K_0}{U_1^-(r_b)})^{\ga}(\fr{\ga+1}{2\ga}U_1^-(r_b)-\fr{\ga-1}{2\ga}\dfr{K_0}{U_1^-(r_b)})
\end{align}
and
\begin{align}\no
\fr{A_0^+}{A_0^-}=\fr{x^\ga}{\ga+1}(-(\ga-1)+\fr{4\ga}{(\ga+1)x-(\ga-1)})=:T_1(x),\\\nonumber
\frac{p_0^+(r_b)}{p_0^-(r_b)}= \fr{x^\ga}{\ga+1}(-(\ga-1)+\fr{4\ga}{(\ga+1)x-(\ga-1)})=: T_2(x),
\end{align}
where $x=\fr{U_1^+(r_b)}{U_1^-(r_b)}\in(\fr{\ga-1}{\ga+1},\fr{\ga+1}{\ga-1})$. It is easy to see that both $T_1(x)$ and $T_2(x)$ monotonically decrease as $x\in(\fr{\ga-1}{\ga+1}, \fr{\ga+1}{\ga-1})$ increases. And $A_0^+>A_0^-, p_0^+(r_b)>p_0^-(r_b)$ if and only if $U_1^-(r_b)>U_1^+(r_b)$.

Next, we consider the relationship between the shock position $r=r_b$ and the exit pressure $p(r_1)=p_{ex}$. In the following, $p_{ex}$ is regarded as a function of $r_b$.

On the exit $r=r_1$,
\be\no
\begin{cases}
p_{ex}=A_0^+(\rho^+(r_1))^{\ga},\\
B_0=\dfr12(\dfr{\ka_1^2}{(\rho^+(r_1))^2r_1^2}+\dfr{\ka_2^2}{r_1^2})+\dfr\ga{\ga-1}A_0^+(\rho^+(r_1))^{\ga-1},
\end{cases}
\ee
which implies,
\be\no
\begin{cases}
\dfr{dp_{ex}}{dr_b}=\dfr{dA_0^+}{dr_b}(\rho^+(r_1))^{\ga}+A_0^+\ga(\rho^+(r_1))^{\ga-1}\dfr{d\rho^+(r_1)}{dr_b},\\
\dfr{d\rho^+(r_1)}{dr_b}=\dfr{\ga (\rho^+(r_1))^{\ga}}{(\ga-1)((U_1^+(r_1))^2-(c^+(r_1))^2)}\dfr{dA_0^+}{dr_b}.
\end{cases}
\ee
Thus,
\be\no
\dfr{dp_{ex}}{dr_b}=(\rho^+(r_1))^{\ga}\dfr{(\ga-1)(M_1^+(r_1))^2+1}{(\ga-1)((M_1^+(r_1))^2-1)}\dfr{dA_0^+}{dr_b}
\ee
with $M_i^\pm(r)=\dfr{U_1^\pm(r)}{c^\pm(r)}$, $i=1,2$.
Note that $(\rho^+(r_1))^{\ga}\dfr{(\ga-1)(M_1^+(r_1))^2+1}{(\ga-1)((M_1^+(r_1))^2-1)}<0$ due to $M_1^+(r_1)<1$.

On the  shock $r=r_b$, one has that
\be\no
B_0=\dfr12(\dfr{\ka_1^2}{(\rho^+(r_b))^2 r_b^2}+\dfr{\ka_2^2}{r_b^2})+\dfr\ga{\ga-1}A_0^+(\rho^+(r_b))^{\ga-1},
\ee
which implies
\be\no
\dfr{dA_0^+}{dr_b}&=&\dfr{\ga-1}{\ga (\rho^+(r_b))^{\ga-1}}\bigg(\dfr{(U_1^+(r_b))^2-(c^+(r_b))^2}{\rho^+(r_b)}\dfr{d\rho^+(r_b)}{dr_b}+\dfr{(U_1^+(r_b))^2+(U_2^+(r_b))^2}{r_b}\bigg)\\\nonumber
&=&(\ga-1)A_0^+\bigg(\dfr{(M_1^+(r_b))^2-1}{\rho^+(r_b)}\dfr{d\rho^+(r_b)}{dr_b}+\dfr{(M_1^+(r_b))^2+(M_2^+(r_b))^2}{r_b}\bigg).
\ee
On the other hand, \eqref{rhop} implies that
\be\no
\dfr{d\rho^+(r_b)}{dr_b}=-\dfr{\ka_1^2}{K_0r_b^2(\rho^-(r_b))^2}\dfr{d\rho^-(r_b)}{dr_b}-\dfr{2\ka_1^2}{K_0r_b^3\rho^-(r_b)}-\dfr{\rho^+(r_b)}{K_0}\dfr{dK_0}{dr_b}.
\ee
By $B_0=\dfr12(\dfr{\ka_1^2}{(\rho^-(r_b))^2 r_b^2}+\dfr{\ka_2^2}{r_b^2})+\dfr\ga{\ga-1}A_0^-(\rho^-(r_b))^{\ga-1}$, one has
\be\no
\dfr{d\rho^-(r_b)}{dr_b}=-\dfr{\fr{\ka_1^2}{\rho^-(r_b)}+\rho^-(r_b)\ka_2^2}{r_b^3((U_1^-(r_b))^2-(c^-(r_b))^2)}.
\ee
Substituting this into the last formula yields
\be\no
\dfr{d\rho^+(r_b)}{dr_b}=\dfr{\rho^+(r_b)}{r_b}\dfr{2+(M_2^-(r_b))^2-(M_1^-(r_b))^2}{(M_1^-(r_b))^2-1}-\dfr{\rho^+(r_b)}{K_0}\dfr{dK_0}{dr_b}.
\ee

Thus,
\begin{align*}
\dfr1{(\ga-1)A_0^+}\dfr{dA_0^+}{dr_b}=&\dfr{(M_1^+(r_b))^2-1}{r_b((M_1^-(r_b))^2-1)}\big(2+(M_2^-(r_b))^2-(M_1^-(r_b))^2\big)\\
&+\dfr{(M_1^+(r_b))^2+(M_2^+(r_b))^2}{r_b}-\dfr{(M_1^+(r_b))^2-1}{K_0}\dfr{dK_0}{dr_b}\\
=&\dfr{1}{r_b\big((M_1^-(r_b))^2-1\big)}\bigg((M_1^+(r_b))^2+(M_1^-(r_b))^2-2\\
&+(M_1^+(r_b))^2(M_2^-(r_b))^2+(M_1^-(r_b))^2(M_2^+(r_b))^2-(M_2^+(r_b))^2-(M_2^-(r_b))^2\bigg)\\
&+\dfr{1-(M_1^+(r_b))^2}{K_0}\dfr{dK_0}{dr_b}.
\end{align*}

We claim that
\be\label{ineq1}
&&(M_1^+(r_b))^2+(M_1^-(r_b))^2-2>0, \\\label{ineq2} &&(M_1^+(r_b))^2(M_2^-(r_b))^2+(M_1^-(r_b))^2(M_2^+(r_b))^2-(M_2^+(r_b))^2-(M_2^-(r_b))^2\geq0,\\\label{ineq3}
&&\dfr{dK_0}{dr_b}\geq0.
\ee
This claim can be verified as follows. Since $K_0=\fr{2(\gamma-1)}{\ga+1}B-\fr{\ga-1}{\ga+1}(U_2^\pm(r_b))^2$, direct calculations show that
\begin{align*}
\dfr{d K_0}{d r_b}=\dfr{2(\ga-1)}{\ga+1}\dfr{(U_2^-(r_b))^2}{r_b}\geq0,
\end{align*}
and
\be\no
\fr1{(M_1^\pm(r_b))^2}=\fr{\ga+1}2 \fr{K_0}{(U_1^\pm(r_b))^2}-\fr{\ga-1}2.
\ee
On the other hand, the inequality \eqref{ineq1} is equivalent to
\be\label{ineq11}
\dfr{2-(M_1^-(r_b))^2}{(M_1^+(r_b))^2}<1.
\ee
Simple computations yield
\be\no
2-(M_1^-(r_b))^2=\dfr{(\ga+1)\fr{K_0}{(U_1^-(r_b))^2}-\ga}{\fr{\ga+1}2 \fr{K_0}{(U_1^-(r_b))^2}-\fr{\ga-1}2},\\\no
\dfr1{(M_1^+(r_b))^2}=\dfr{\ga+1}2 \dfr{(U_1^-(r_b))^2}{K_0}-\dfr{\ga-1}2.
\ee
Note that $\fr{U_1^+(r_b)}{U_1^-(r_b)}>1>\fr{\ga-1}{\ga+1}$, thus $\dfr{\ga+1}2 \dfr{K_0}{(U_1^-(r_b))^2}-\dfr{\ga-1}2>0$. Thus \eqref{ineq11} is equivalent to
\be\no
\bigg((\ga+1)\dfr{K_0}{(U_1^-(r_b))^2}-\ga\bigg)\bigg(\dfr{\ga+1}2 \dfr{(U_1^-(r_b))^2}{K_0}-\dfr{\ga-1}2\bigg)<\dfr{\ga+1}2 \dfr{K_0}{(U_1^-(r_b))^2}-\dfr{\ga-1}2.
\ee
That is,
\be\no
\dfr{K_0}{(U_1^-(r_b))^2}+\dfr{(U_1^-(r_b))^2}{K_0}>2,
\ee
which is always valid. Thus \eqref{ineq1} holds. To verify \eqref{ineq2}, one needs only to consider $\ka_2\neq 0$ since $\kappa_2=0$ is trivial. Then the inequality \eqref{ineq2} is equivalent to\\
\be\no
\dfr{(M_1^-(r_b))^2-1}{1-(M_1^+(r_b))^2}>\dfr{(M_2^-(r_b))^2}{(M_2^+(r_b))^2}=\dfr{(c^+(r_b))^2}{(c^-(r_b))^2},
\ee
i.e.,
\be\label{ineq21}
(U_1^-(r_b))^2-(c^-(r_b))^2>(c^+(r_b))^2-(U_1^+(r_b))^2.
\ee
Since $K_0=\fr{2}{\ga+1}(c^\pm(r_b))^2+\fr{\ga-1}{\ga+1}(U_1^\pm(r_b))^2$, so \eqref{ineq21} is equivalent to
\be\no
(U_1^-(r_b))^2-\dfr{\ga+1}2 K_0+\dfr{\ga-1}2 (U_1^-(r_b))^2>\dfr{\ga+1}2 K_0-\dfr{\ga-1}2 (U_1^+(r_b))^2-(U_1^+(r_b))^2
\ee
That is,
\be\no
2 K_0<(U_1^+(r_b))^2+(U_1^-(r_b))^2=(U_1^-(r_b))^2+\dfr{K_0^2}{(U_1^-(r_b))^2},
\ee
which is always valid. Thus \eqref{ineq2} holds true.



%


In summary, we have shown that $\dfr{dA_0^+}{dr_b}>0$ and thus, $\dfr{dp_{ex}}{dr_b}<0$, which implies that the shock position $r=r_b$ increases as the exit pressure $p_{ex}$ decreases. So the following result holds.
\bp\label{theorem1}
{\it For any given incoming flow $U_1^-(r_0)$, $U_2^-(r_0)$, $\rho^-(r_0)$, $p^-(r_0)$ satisfying $U_1^-(r_0)>c^-(r_0)$, there exist two constants $p_1<p_0$ such that, for any $p_{ex}\in(p_1,p_0)$, there exists a unique solution to Problem III with a shock located at $r\equiv r_b$ satisfying $U_1^-(r)>c^-(r)$ for $r\in(r_0,r_b)$, $U_1^+(r)<c^+(r)$ for $r\in(r_b,r_1)$. In addition, the shock position $r=r_b$ increases as the exit pressure $p_{ex}$ decreases. Furthermore, the shock position $r_b$ approaches to $r_1$ if $p_{ex}$ goes to $p_1$ and $r_b$ tends to $r_0$ if $p_{ex}$ goes to $p_0$.
}\ep

It remains to clarify the flow state just behind the shock. That is, to determine the sign of
\be\no
(U_{1}^+)^2+(U_{2}^+)^2 -c^2(\rho^+,A_0^+)
\ee
at $r=r_b$. It follows from \eqref{Omgp}-\eqref{rhop} that
\begin{eqnarray}\nonumber
&&(U_{1}^+)^2+ (U_{2}^+)^2 -c^2(\rho^+,A_0^+) =U_{1}^2-\frac{\gamma+1}{2}(K_0-\frac{\gamma-1}{\gamma+1}U_{1}^2) + U_{2}^2\\\no
&&=\frac{\gamma+1}{2}(U_{1}^2-K_0)+U_{2}^2=\frac{\gamma+1}{2}\frac{K_0}{(U_{1}^-)^2}(K_0-(U_{1}^-)^2)+ U_{2}^2\\\no
&&=\frac{K_0}{(U_{1}^-)^2}\left(c^2(\rho^-, A_0^-)-(U_{1}^-)^2+\frac{(U_{1}^- U_{2}^-)^2}{K_0}\right)\\\no
&&=\frac{K_0}{(M_{1}^-)^2}\bigg(1-(M_{1}^-)^2+\frac{(M_{1}^-)^2(M_{2}^-)^2}{\frac{2}{\gamma+1}+\frac{\gamma-1}{\gamma+1}(M_{1}^-)^2}\bigg)\\\no
&&=\frac{K_0}{(M_{1}^-)^2(\frac{2}{\gamma+1}+\frac{\gamma-1}{\gamma+1}(M_{1}^-)^2)}\bigg((1-(M_{1}^-)^2)(\frac{2}{\gamma+1}+\frac{\gamma-1}{\gamma+1}(M_{1}^-)^2)+(M_{1}^-)^2(M_{2}^-)^2\bigg)
\end{eqnarray}
Since $B_0=\fr12((U_1^-)^2+(U_2^-)^2)+\fr{(c^-)^2}{\ga-1}$, thus
\be\no
(M_2^-)^2=\fr{(U_2^-)^2}{2B_0-(U_2^-)^2}((M_1^-)^2+\fr2{\ga-1}).
\ee

Assume that $\kappa_2\neq0$. Define
\begin{eqnarray}\nonumber
&&x(r)=(M_{1}^-(r))^2,\ \ a(r)=\fr{\ka_2^2}{2r^2B_0-\ka_2^2}>0,\\\nonumber
&& g_{a}(x)=(1-x)(\f{2}{\ga+1}+\f{\ga-1}{\ga+1}x)+ax(x+\fr2{\ga-1})\\\no
&&=\fr2{\ga+1}+(\fr{2a}{\ga-1}+\fr{\ga-3}{\ga+1})x+(a-\fr{\ga-1}{\ga+1})x^2.
\end{eqnarray}
The definition of $f_1$ in \eqref{f1} shows $f_1(r)=(1-\frac{\ga+1}{\ga-1}a(r))x(r)$. To clarify different cases, one introduces two positive constants $r'_*>r_*>0$, which are determined as follows. Since
\be\label{ainfty}
a'(r)<0,\ \ \lim_{r\rightarrow +\infty}a(r)=0,
\ee
there exists a unique $r_*=\sqrt{\frac{\ga \kappa_2^2}{(\ga-1)B_0}}$ such that $a(r_*)=\frac{\ga-1}{\ga+1}$ .

Since $a(r)< \frac{\ga-1}{\ga+1}$ for any $r> r_*$ and decreases as $r$ increases, $x(r)>1$ and $x'(r)>0$ for any $r>r_0$, one has that $f_1'(r)>0$ for all $r>r_*$. It follows from \eqref{M1infty} and \eqref{ainfty} that
\be\nonumber
\lim_{r\rightarrow +\infty}f_1(r)=+\infty.
\ee
Since $f_1(r_*)=0$, there exists a unique $r'_*>r_*$ such that $f_1(r'_*)=1$.

To distinguish whether the state behind the shock is subsonic or supersonic, it suffices to check whether $g_{a(r_b)}(x(r_b))>0$ or $<0$. It is easy to find the roots of $g_{a}(x)=0$, which are denoted by
\begin{eqnarray}\nonumber
x_1(r)\equiv-\frac2{\ga-1}<0,\quad x_2(r)=\frac1{1-\frac{\ga+1}{\ga-1}a(r)}.
\end{eqnarray}
First, we determine the sign of $g_{a(r_b)}(x(r_b))$ by the values of $a$ and $f_1$ at $r=r_b$ as follows.
\begin{enumerate}[(1)]
  \item $a(r_b)>\frac{\ga-1}{\ga+1}$. Then $x_1<x_2(r_b)<0$, and $g_{a(r_b)}(x(r_b))\geq g_a(0)=\frac2{\ga+1}>0$ for all $x\geq0$, which implies that the flow behind the shock is supersonic.
  \item $a(r_b)=\frac{\ga-1}{\ga+1}$. Then $g_{a(r_b)}(x(r_b))=\frac2{\ga+1}+\fr{\ga-1}{\ga+1}x(r_b)\geq\frac2{\ga+1}>0$ for all $x(r_b)\geq0$. Thus the flow behind the shock is supersonic also.
  \item $a(r_b)<\frac{\ga-1}{\ga+1}$. Then $x_1=-\frac2{\ga-1}<0<1<x_2(r_b)=\frac1{1-\frac{\ga+1}{\ga-1}a(r_b)}$, and $g_{a(r_b)}(x(r_b))>0$ for $0\leq x(r_b)<x_2(r_b)$, while $g_{a(r_b)}(x(r_b))<0$ for $x(r_b)>x_2(r_b)$. Thus, it follows that
      \begin{enumerate}[(i)]
      \item if $f_1(r_b)<1$, then the state behind the shock is supersonic;
      \item if $f_1(r_b)>1$, then the state behind the shock is subsonic;
      \item if $f_1(r_b)=1$, then the state behind the shock is sonic.
      \end{enumerate}
\end{enumerate}

Since the position $r_b$ of the shock is unknown, one needs to use the properties of $a$, $f_1$ and the monotonicity between the shock position and the exit pressure to determine the sign of $g_{a(r_b)}(x(r_b))$ in the following cases:

{\bf Case a}: $r_0<r_*$ (i.e. $a(r_0)>\fr{\ga-1}{\ga+1}$).  The flow is supersonic (subsonic) behind the shock if the shock occurs at $r_0<r_b<r'_*$ ($r_b>r'_*$), and is sonic behind the shock if $r_b=r'_*$. Since the exit pressure $p_{ex}$ depends monotonically on the shock position $r=r_b$, there exists an interval $(p_1,p_0)$ such that if the exit pressure $p_{ex}\in(p_1,p_0)$, there exists a unique shock solution with a shock at $r_b\in(r_0,r_1)$.

If $r_1\leq r'_*$, one obtains a supersonic-supersonic shock. Note that the flow at downstream can change smoothly from supersonic to subsonic after crossing the shock.

If $r_1> r'_*$, then there exists a $p'_*\in(p_1,p_0)$ corresponding to $r'_*$ such that
\begin{enumerate}[(i)]
  \item if $p_{ex}\in(p_1,p'_*)$, one has a supersonic-subsonic shock located at $r_b\in(r'_*,r_1)$;
  \item if $p_{ex}\in(p'_*,p_0)$, one gets a supersonic-supersonic shock located at $r_b\in(r_0,r'_*)$ and the flow at downstream can change smoothly from supersonic to subsonic after crossing the shock.
  \item if $p_{ex}=p'_*$, one obtains a supersonic-sonic shock with a shock and the sonic circle both located at $r_b=r'_*$.
\end{enumerate}

{\bf Case b}:  $r_*\leq r_0<r'_*$.  The flow is supersonic behind the shock if the shock occurs at $r_0<r_b<r'_*$, and is subsonic (sonic) behind the shock if the shock occurs at $r_b>r'_* (r_b=r'_*)$. Therefore there exists an interval $(p_1,p_0)$ such that if the exit pressure $p_{ex}\in(p_1,p_0)$, there exists a unique shock solution with a shock at $r_b\in(r_0,r_1)$.

  If $r_1\leq r'_*$, then the shock is supersonic-supersonic, the flow at downstream can change smoothly from supersonic to subsonic crossing the shock.

  If $r_1> r'_*$, then there exists a $p'_*\in(p_1,p_0)$ corresponding to $r'_*$ such that
  \begin{enumerate}[(i)]
  \item if $p_{ex}\in(p_1,p'_*)$, one gets a supersonic-subsonic shock with the shock located at $r_b\in(r'_*,r_1)$;
  \item if $p_{ex}\in(p'_*,p_0)$, one has a supersonic-supersonic shock with the shock located at $r_b\in(r_0,r'_*)$, the flow at downstream can change smoothly from supersonic to subsonic crossing the shock;
  \item if $p_{ex}=p'_*$, one obtains a supersonic-sonic shock with the shock and the sonic line both located at $r_b=r'_*$.
  \end{enumerate}

{\bf Case c}: $r_0\geq r'_*>r_*$. Since $f_1(r_0)\geq f_1(r'_*)=1$, $f_1(r)>1$ for any $r>r_0$. The flow is subsonic behind the shock if the shock occurs at $r_b>r_0$. Therefore, there exists an interval $(p_1,p_0)$ such that, for any $p_{ex}\in(p_1,p_0)$, there exists a unique supersonic-subsonic shock solution with the shock located at $r_b\in(r_0,r_1)$.

In summary, we have proved
\bp\label{prop1}
{\it For the solutions in Proposition \ref{theorem1}, the flow state just behind the shock can be divided into the following three cases.
\begin{enumerate}[1.]
  \item If $r_0<r_1\leq r'_*$, then $|{\bf M}^+|^2(r_b)>1$ for any $p_{ex}\in(p_1,p_0)$.
  \item If $r_0<r'_*<r_1$, there exists a $p'_*\in(p_1,p_0)$ such that $|{\bf M}^+|^2(r_b)>1$ for any $p_{ex}\in(p'_*,p_0)$, $|{\bf M}^+|^2(r_b)<1$ for any $p_{ex}\in(p_1,p'_*)$ and $|{\bf M}^+|^2(r_b)=1$ for $p_{ex}=p'_*$.
  \item If $r_0\geq r'_*$, then $|{\bf M}^+|^2(r_b)<1$ for any $p_{ex}\in(p_1,p_0)$.
  \end{enumerate}
}
\ep


Finally one can determine the flow state at the exit by examing the sign of
\be\no
(U_{1}^+)^2+(U_{2}^+)^2 -c^2(\rho^+,A_0^+)
\ee
at $r=r_1$.

Since $\ka_1$, $\ka_2$ and $B_0$ are unchanged across the shock, it follows from the Bernoulli's law and the exit pressure $p_{ex}=A_0^+(\rho^+(r_1))^\ga$ that $\rho^+(r_1)$ satisfies $f_2(\rho^+(r_1))=0$, where $f_2$ is defined in \eqref{f2}.

Since $f'_2(\rho)=\dfr{\ga p_{ex}}{\ga-1}-2(B_0-\dfr12\dfr{\ka_2^2}{r_1^2})\rho$, $f_2(\rho)$ increases on $(0,\rho_\sharp)$ and decreases on $(\rho_\sharp,\infty)$ as $\rho$ increases, where
\be\nonumber
\rho_\sharp=\dfr{\ga p_{ex}}{2(\ga-1)(B_0-\fr12\fr{\ka_2^2}{r_1^2})}.
\ee
Note that $f_2(0)>0$, thus $\rho<\rho^+(r_1)$ if and only if $f_2(\rho)>0$ and $\rho>\rho^+(r_1)$ if and only if $f_2(\rho)<0$.

The flow is supersonic at $r=r_1$ if and only if $\dfr{2(\ga-1)}{\ga+1}B_0>\dfr{\ga p_{ex}}{\rho^+(r_1)}$, that is,
\be\label{rhosharp}
\rho^+(r_1)> \dfrac{\ga(\ga+1) p_{ex}}{2(\ga-1)B_0}=:\rho^{\sharp},
\ee
which is equivalent to $f_2(\rho^\sharp)>0$.
And the flow is subsonic at $r=r_1$ if and only if $\rho^+(r_1)<\rho^\sharp$, which is equivalent to $f_2(\rho^\sharp)<0$.

In summary, we have shown the following proposition.
\bp\label{prop2}
{\it For the solutions given in Proposition \ref{theorem1}, the flow state at the exit can be classified as:
\begin{enumerate}[1.]
  \item If $f_2(\rho^\sharp)>0$, then $|{\bf M}^+|^2(r_1)>1$.
  \item If $f_2(\rho^\sharp)<0$, then $|{\bf M}^+|^2(r_1)<1$.
  \item If $f_2(\rho^\sharp)=0$, then $|{\bf M}^+|^2(r_1)=1$.
\end{enumerate}
}\ep
By Propositions \ref{theorem1}, \ref{prop1} and \ref{prop2}, the proof of Theorem \ref{ProbIII} is completed.

\subsection{Transonic shock flows moving from the outer to inner circle}\noindent

Finally, we consider the transonic shock wave patterns when a supersonic flow moves from the outer to the inner circle.

{\bf Problem IV.} Construct smooth functions $(U_{1}^{-},U_{2}^{-}, \rho^{-}, p^-, A^-)$ and $(U_{1}^{+},U_{2}^{+}, \rho^{+}, p^+, A^+)$, which satisfy the system \eqref{radial-euler} in $\Omega^-=\{(x_1,x_2): r_b<r<r_1\}$ and $\Omega^+=\{(x_1,x_2): r_0<r<r_b\}$ respectively, supplemented with the boundary conditions
\begin{eqnarray}\label{IV-bd1}
&&\rho(r_1)=\rho_{0}^-,\ U_{1}(r_1)=U_{10}^-<0, \ U_{2}(r_1)=U_{20}^-\neq 0,\ A(r_1)=A_{0}^-,\\\label{IV-bd2}
&&p(r_0)= p_{ex}.
\end{eqnarray}
The discontinuity occurs at $r=r_b\in (r_0,r_1)$ which is unknown, and across the discontinuity $r=r_b$, the Rankine-Hugoniot condition \eqref{rankine-hugoniotII} holds and the entropy increases.

In this case, $\kappa_1= r_1 \rho_0 U_{10}^-$ and $\kappa_2=r_1 U_{20}^-\neq 0$. Recall that
\begin{eqnarray}\nonumber
f_1(r)&=&\bigg(1-\frac{(\gamma-1)\kappa_2^2}{(\gamma+1)(2r^2 B_0-\kappa_2^2)}\bigg)(M_1^-)^2(r),\\\nonumber
\rho^{\sharp}&=&\dfrac{\ga(\ga+1) p_{ex}}{2(\ga-1)B_0}.
\end{eqnarray}
and modify $f_2$ as
\begin{eqnarray}\nonumber
f_2(\rho)&=&\dfr{\ga p_{ex}}{\ga-1}\rho-(B_0-\dfr12\dfr{\ka_2^2}{r_0^2})\rho^2+\dfr{\ka_1^2}{2r_0^2},
\end{eqnarray}
Denote the limit circle defined in {\bf Problem II} with boundary condition \eqref{IV-bd1} at $r=r_1$ by $r=r^{\sharp,-}$, which is $(U_1^-)^2(r^{\sharp,-})= (c^-)^2(r^{\sharp,-})$. Therefore $r^{\sharp,-}$ is the unique root of the following equation in $\mathbb{R}^+$
\be\label{limit-eq}
\fr{2(\ga-1)B_0}{\ga+1}x^2=(\ga A_0^-)^{\fr2{\ga+1}}\kappa_1^{\fr{2(\ga-1)}{\ga+1}}x^{\fr4{\ga+1}}+\fr{\ga-1}{\ga+1}\kappa_2^2.
\ee

Denote the limit circle defined in {\bf Problem II} with boundary data $(U_1^+(r_b), U_1^+(r_b), \rho^+(r_b), A_0^+)$ at $r=r_b$ by $r=r^{\sharp,+}$, then $r^{\sharp,+}$ is the unique root to \eqref{limit-eq} with $A_0^-$ replaced by $A_0^+$. Since $A_0^+>A_0^-$, it is easy to verify that $r^{\sharp,-}<r^{\sharp,+}$.

First, to compare $r_0$ and $r^{\sharp,+}$, one defines
\be\label{rhoshatpsharp}
\rho^{\sharp\sharp}=\dfr{\ga p_{ex}}{K(r_0)}\geq \rho^\sharp.
\ee
The flow is supersonic in $r$-direction at $r=r_0$ if and only if $K(r_0)>\dfr{\ga p_{ex}}{\rho^+(r_0)}$, that is,
\be\label{rhosharpsharp}
\rho^+(r_0)> \rho^{\sharp\sharp},
\ee
which is equivalent to $f_2(\rho^{\sharp\sharp})>0$.
And the flow is subsonic in $r$-direction at $r=r_0$ if and only if $\rho^+(r_0)<\rho^{\sharp\sharp}$, which is equivalent to $f_2(\rho^{\sharp\sharp})<0$. Therefore, $r_0\geq r^{\sharp,+}$ if and only if $f_2(\rho^{\sharp\sharp})\leq0$, which is independent of the location of the shock $r_b$.

For fixed $r_0\in (r^{\sharp,+}, r_1)$, the monotonicity between the shock position $r=r_b$ and the exit pressure $p(r_0)=p_{ex}$ is also valid in this case, that is, $\dfr{dp_{ex}}{dr_b}<0$. Then there exist two positive constants $p_0=p_{ex}(r_0)>p_1=p_{ex}(r_1)$ such that if the exit pressure $p_{ex}$ belongs to $(p_1,p_0)$, there exists a piecewise smooth solution to {\bf Problem IV} with a shock $r=r_b\in (r_0,r_1)$.

To determine the flow patterns of the shock $r=r_b$, as in {\bf Problem III}, one needs to introduce some parameters. Since $f_1(r)$ is indeed well-defined in $(r^{\sharp,-},\infty)$, and $(M_1^-)^2(r)>1$ and $\frac{d}{dr}(M_1^-)^2(r)>0$ for any $r>r^{\sharp,-}$, one has $f_1'(r)>0$ for all $r>r_*=\sqrt{\frac{\gamma}{(\gamma-1) B_0}}\kappa_2$. Due to $f_1(r_*)=0$ and $\lim_{r\rightarrow +\infty}f_1(r)=+\infty$, there exists a unique $r'_*>r_*$ such that $f_1(r'_*)=1$.

%

Then we can obtain the following theorem by using Propositions \ref{prop1} and \ref{prop2}.

\bt\label{ProbIV}
{\it For a given incoming flow, which is supersonic in the
$r$-direction with nonzero angular velocity at the entrance $r=r_1$ and exit pressure $p_{ex}$ at $r=r_0$, then there exist two positive constants $p_0<p_1$ such that if $f_2(\rho^{\sharp\sharp})>0$ or $f_2(\rho^{\sharp\sharp})\leq0$ but $p_{ex}\notin(p_1,p_0)$, there does not exist any piecewise smooth solution to the Problem IV in $\Omega$; if $f_2(\rho^{\sharp\sharp})\leq0$ and $p_{ex}\in(p_1,p_0)$, there exists a unique piecewise smooth weak solution to Problem IV in $\Omega$ with a shock located at $r= r_b\in(r_0,r_1)$ with the following properties.
\begin{enumerate}[(A)]
  \item The shock position $r=r_b$ increases as the exit pressure $p_{ex}$ decreases. In addition, the shock position $r_b$ tends to $r_1$ if $p_{ex}$ goes to $p_1$ and $r_b$ approaches to $r_0$ if $p_{ex}$ goes to $p_0$.
  \item The flow patterns in $\Omega^+$ can be classified in terms of the boundary conditions as follows.
  \begin{enumerate}[1.]
  \item If $r_1>r_0\geq r'_*$, then
  \begin{description}
  \item{Subcase 1.1.} for $f_2(\rho^\sharp)<0$, it is a supersonic-subsonic shock and the flow is uniformly subsonic in $\Omega^+$.
  \item{Subcase 1.2.} for $f_2(\rho^\sharp)>0$, it is a supersonic-subsonic shock and the flow changes smoothly from subsonic to supersonic in $\Omega^+$.
  \item{Subcase 1.3.} for $f_2(\rho^\sharp)=0$, it is a supersonic-subsonic shock and the flow is subsonic in $\Omega^+$ but degenerates to sonic state at the exit.
  \end{description}
  \item If $r_1>r'_*> r_0$, there exists a $p'_*\in(p_1,p_0)$ such that,
  \begin{description}
  \item{Subcase 2.1.} for $p_{ex}\in(p_1,p'_*)$ with $f_2(\rho^\sharp)<0$, it is a supersonic-subsonic shock and the flow is uniformly subsonic in $\Omega^+$.
  \item{Subcase 2.2.} for $p_{ex}\in(p_1,p'_*)$ with $f_2(\rho^\sharp)>0$, it is a supersonic-subsonic shock and the flow changes smoothly from subsonic to supersonic in $\Omega^+$.
  \item{Subcase 2.3.} for $p_{ex}\in(p_1,p'_*)$ with $f_2(\rho^\sharp)=0$, it is a supersonic-subsonic shock and the flow is subsonic in $\Omega^+$ but degenerates to sonic state at the exit.
  \item{Subcase 2.4.} for $p_{ex}\in(p'_*,p_0)$, there exists a supersonic-supersonic shock and the flow is uniformly supersonic in $\Omega^+$.
  \item{Subcase 2.5.} for $p_{ex}=p'_*$, there is a supersonic-sonic shock and the flow is supersonic in $\Omega^+$ but degenerates to the sonic state at the shock position.
  \end{description}
  \item If $r_1 \leq r'_*$, then, for any $p_{ex}\in(p_1,p_0)$, there exists a supersonic-supersonic shock and the flow is uniformly supersonic in $\Omega^+$.
\end{enumerate}
\end{enumerate}
}
\et

\br
{\it Since $M^2_1(r'_*)=\fr1{1-(\ga+1)a(r'_*)/(\ga-1)}\geq 1=M_1^2(r^{\sharp,+})$, thus $r'_*\geq r^{\sharp,+}$  and ``=" is reached if and only if $\kappa_2=0$, which is equivalent to $\rho^{\sharp\sharp}=\rho^\sharp$. Thus, only the subcase 1.1 and 1.3 in Theorem \ref{ProbIV} can occur when $\kappa_2=0$ in {\bf Problem IV.}
}
\er

{\bf Acknowledgement.}  Part of this work was done when Weng visited The Institute of Mathematical Sciences of The Chinese University of Hong Kong. He is grateful to the institute for providing nice research environment. Weng is partially supported by National Natural Science Foundation of China 11701431, 11971307, 12071359, the grant of Project of Thousand Youth Talents (No. 212100004). Xin is supported in part by the Zheng Ge Ru  Foundation, Hong Kong RGC Earmarked Research Grants CUHK-14305315, CUHK-14300917, CUHK-14302819 and CUHK-14302917, and by Guangdong Basic and Applied Basic Research Fundation 2020B1515310002.


\end{document}